%  TEXSTUFF +++++++++++++++++++++++++++++++++++++++++++++++++++++++++++++
\input amssym.def
\input epsf

\let \blskip = \baselineskip
\parskip=1.2ex plus .2ex minus .1ex

\tabskip 20pt
\tolerance = 1000
\pretolerance = 50
\newcount\itemnum
\itemnum = 0
\overfullrule = 0pt

% HEADER DEFNS ++++++++++++++++++++++++++++++++++++++++++++++++++++
\def\title#1{\bigskip\centerline{\bigbigbf#1}}
\def\author#1{\bigskip\centerline{\bf #1}\smallskip}
\def\address#1{\centerline{\it#1}}
\def\abstract#1{\vskip1truecm{\narrower\noindent{\bf Abstract.} #1\bigskip}}

%  GENERAL MACROS    +++++++++++++++++++++++++++++++++++++++++++++++++++++
\def\sp{\bigskip}
\def\nosp{\vskip -\the\blskip plus 1pt minus 1pt}
\def\upsp{\nosp\medskip}
\def\br{\hfil\break} 
\def\ti{\br \hglue \the \parindent}

\def\ce#1{\LP\centerline{#1}}

\def\skipit#1{}
\def\mdag{\raise 3pt\hbox{\dag}}

\def\XP{\par\noindent\hang}
\def\LP{\par\noindent}
\def\BP[#1]{\par\item{[#1]}}
\def\SH#1{\sp\vskip\parskip\leftline{\bigbf #1}\nobreak}

\def\TH#1{\sp\XP{\bf THEOREM\ \shead#1}}
\def\LM#1{\sp\XP{\bf LEMMA\ \shead#1}}

\def\PF{\LP{\bf Proof:\ }}
\def\NX{\advance\itemnum by 1 \sp\LP {\bf \shead \the\itemnum.\ }}
\def\qed{\null\nobreak\hfill\hbox{${\vrule width 5pt height 6pt}$}\par\sp}

\def\cart{\>\hbox{${\vcenter{\vbox{
    \hrule height 0.4pt\hbox{\vrule width 0.4pt height 4.5pt
    \kern4pt\vrule width 0.4pt}\hrule height 0.4pt}}}$}\>}
\def\bxmu{\>\hbox{${\vcenter{\vbox {
    \hrule height 0.4pt\hbox{\vrule width 0.4pt height 4pt
    \hskip -1.3pt\lower 1.8pt\hbox{$\times$}\negthinspace\vrule width 0.4pt}
    \hrule height 0.4pt}}}$}\>}

\def\lin#1{\hbox to #1true in{\hrulefill}}

% STYLE COMMANDS +++++++++++++++++++++++++++++++++++++++++++++++++++++++

\def\foot#1{\raise 6pt \hbox{#1} \kern -3pt}
% Printed \today.}

% FIGURES AND TABLES +++++++++++++++++++++++++++++++++++++++++++++++++++++++++++
\def\fig #1 #2 #3 #4 #5 {\sp \ce{ {\epsfbox[#1 #2 #3 #4]{figs/#5.ps}} }}
%left x, bottom y, right x, top y, in portrait orientation

%\def\gpic#1{#1 \sp\ce{\box\graph} \medskip} %gpic picture, centered with space
%for troff tables: \s\ccol = centered col, \s\vcol = rule, # to end preamble

% JOURNAL ABBREVS  +++++++++++++++++++++++++++++++++++++++++++++++++++++

\def\JGT{{\it J.\ Graph Theory}}

\def\DM{{\it Discrete Math.{}}}

\def\SIAD{{\it SIAM\ J.\ Algeb.\ Disc.\ Meth.{}}}

\def\SIDM{{\it SIAM\ J.\ Discr.\ Math.{}}}

\def\CNum{{\it Congr.\ Numer.{}}}

% SPECIAL CHARACTERS ++++++++++++++++++++++++++++++++++++++++++++++++++++

    \def\bs{\backslash}

%\def\square{\hbox{${\vcenter{\vbox{
%    \hrule height 0.3pt\hbox{\vrule width 0.3pt height 4.5pt
%    \kern 3.8pt\vrule width 0.3pt}\hrule height 0.3pt}}}$}}

  %{roman {I back 20 N}}
  %{roman {C back 50 C}}
  %{roman {C back 50 Q}}
  %{\rm \kern -8pt Z} }  %\bold Z}   %{roman {Z back 50 Z}}
  %{roman {I back 20 R}}

% BINARY RELATIONS ++++++++++++++++++++++++++++++++++++++++++++++++++++++
\def\nul{\hbox{\O}}    
		
\def\esub{\subseteq}

		%also \vee,\wedge

    %Note also uppercase doubles line, plus [long][left][right]arrow
    %and [up][down]arrow, plus ne,sw,se,nw+arrow,
    %plus [long]mapsto, hook[left][right]arrow, and harpoons

% DELIMITERS ++++++++++++++++++++++++++++++++++++++++++++++++++++++++++++++
\def\({\left(}	\def\){\right)}

% MACROS WITH ARGUMENTS +++++++++++++++++++++++++++++++++++++++++++++++

\def\CH#1#2{{{#1}\choose{#2}}}

\def\FL#1{\left\lfloor{#1}\right\rfloor}
\def\CL#1{\left\lceil{#1}\right\rceil}

\def\VEC#1#2#3{#1_{#2},\ldots,#1_{#3}}

\def\LIST#1#2{#1,\ldots,#2}

\def\st{\colon\;} %"such that"

\def\SET#1:#2{\{#1\colon\;#2\}}

% MODES ++++++++++++++++++++++++++++++++++++++++++++++++++++++++++++++++
\def\B#1{{\bf #1}}		
    %cardinality
\def\ov#1{\overline{#1}}

   \def\Db{\ov{D}} \def\Eb{\ov{E}}
\def\Fb{\ov{F}}   \def\Ib{\ov{I}} 
  \def\Mb{\ov{M}}

% WORDS ++++++++++++++++++++++++++++++++++++++++++++++++++++

% GRAPHS +++++++++++++++++++++++++++++++++++++++++++++++++++++++++++

% POSETS ++++++++++++++++++++++++++++++++++++++++++++++++++++++++++

% MATROIDS +++++++++++++++++++++++++++++++++++++++++++++++++++++++++++

% ENUMERATION ++++++++++++++++++++++++++++++++++++++++++++++++++++++

% PAGE SETUP +++++++++++++++++++++++++++++++++++++++++++++++++++++++++
\magnification=\magstep1
\vsize=9.0 true in
\hsize=6.5 true in
%%\hoffset= -.3 true in
%%\voffset= -.1 true in
\headline={\hfil\ifnum\pageno=1\else\folio\fi\hfil}
%\footline={\hfil-- \folio\ --\hfil}
\footline={\hfil\ifnum\pageno=1\folio\else\fi\hfil}

\parindent=20pt
\baselineskip=12pt
\parskip=.5ex  %1.75

% MACRO STUFF (goes to another file?) ++++++++++++++++++++++++++
\def\shead{ }

\font\bigbf = cmb10 scaled \magstep1

\font\bigbigbf = cmb10 scaled \magstep2

% TITLE PAGE FORMAT ++++++++++++++++++++++++++++++++++++++++++++++++
%\def\titpage{\vbox to \vsize
%     {\vfil \title \vfil \author \vfil \vfil \setabs \vfil \extra}  }

%\input eplain
% gpic macros for use with "gpic -t" as preprocessor for tex,
% by Douglas B. West
\def\gpic#1{#1
     \bigskip\par\noindent{\centerline{\box\graph}}
     \medskip} %gpic picture, centered with space

% HEADER SECTION ++++++++++++++++++++++++++++++++++++++++++++++++++++
\title{INTERSECTION REPRESENTATION OF DIGRAPHS}
\title{IN TREES WITH FEW LEAVES}
\author{In-Jen Lin\foot{1}}
\address{National Ocean University, Taipei, Taiwan, ijlin@tiger.cs.nthu.edu.tw}
\author{Malay K. Sen\foot{2}}
\address{North Bengal University, Darjeeling, India}
\author{Douglas B. West\foot{3}}
\address{University of Illinois, Urbana, IL 61801-2975, west@math.uiuc.edu}
\vfootnote{}{\br
  \foot{1}Supported by a University of Illinois Research Board Grant.\br
  \foot{2}Supported by the Council for International Exchange of Scholars.\br
  \foot{3}Supported by NSA/MSP Grant MDA904-93-H-3040.\br
   Running head: {REPRESENTATION OF DIGRAPHS} \br
   Keywords: leafage, digraph, intersection representation, Ferrers dimension\br
   Written September 1992, revised July 1997.
}
\abstract{The {\it leafage} of a digraph is the minimum number of leaves in a
  host tree in which it has a subtree intersection representation.  We discuss
  bounds on the leafage in terms of other parameters (including Ferrers
  dimension), obtaining a string of sharp inequalities.}

% ABBREVIATIONS ++++++++++++++++++++++++++++++++++++++++++++++++++++++++
\def\Ib{\bar I}
\def\Mb{\bar M}
\def\int{intersection}
\def\intg{intersecting}
\def\rep{representation}
\def\dig{digraph}
\def\idig{interval digraph}
\def\fdig{Ferrers digraph}
\def\fd{Ferrers dimension}

\def\charzn{characterization}

% DOCUMENT ++++++++++++++++++++++++++++++++++++++++++++++++++++++++++
\SH
{1. INTRODUCTION}
An {\it \int\ \rep} of a \dig\ $D$ assigns an ordered pair $(S_v,T_v)$
to each vertex $v\in V(D)$ such that $uv \in E(D)$ if and only if
$S_u\cap T_v \ne \nul$.  We call $S_v$ and $T_v$ the {\it source set} and 
{\it sink set} of $v$.  This model was first described by Beineke and Zamfirescu
[1] under the name {\it connection digraph}.  An essentially equivalent model
in terms of bipartite graphs was introduced by Harary, Kabell, and McMorris
[7].

When each set in an \int\ \rep\ is a subtree of a fixed host tree, we have a
{\it subtree \rep}.  Every $n$-vertex \dig\ has a subtree \rep\ in a star
with $n$ leaves.  Not every digraph has a subtree \rep\ in a path; those
that do are the interval digraphs, which are characterized in [15,16].
We define the {\it leafage} $l(D)$ of a \dig\ $D$ to be the minimum number
of leaves in a host tree in which $D$ has a subtree \rep.  Thus leafage is a
measure of distance from an \idig, and the subtree \rep s in stars show that
$l(D)\le n(D)$.  An analogous parameter for chordal (undirected) graphs is
studied in [8].  Further results about adjacency matrices of \idig s appear
in [9,10,15,16,17].

We obtain lower bounds on leafage using the idea of \fd.  The {\it successors}
of a vertex $v$ are $\{u\in V(D)\st vu\in E(D)\}$; the {\it predecessors} are
$\{u\in V(D)\st uv\in E(D)\}$.  A \dig\ is a {\it \fdig} [14] if its successor
sets are linearly ordered by inclusion, which is equivalent to the adjacency
matrix $A(D)$ having no 2 by 2 permutation submatrix.  Viewing a \dig\
as a relation $D\esub V(D)\times V(D)$, the {\it \fd} $f(D)$ of $D$ is the
minimum number of \fdig s on $V(D)$ whose \int\ is $D$ (introduced in [2]).
Since the complement in $V\times V$ of a \fdig\ is also a \fdig, this also
equals the minimum number of \fdig s whose union is $\Db$.

Interval digraphs all have \fd\ at most 2; a \dig\ $D$ is an \idig\ if and only
if $\Db$ is the union of two disjoint \fdig s [15].  This generalizes to a
lower bound on $l(D)$ using Ferrers digraphs.  Let $f^*(D)$ denote the minimum
number of pairwise disjoint \fdig s whose union is $\Db$.  These are \fdig s
whose intersection is $D$ and whose pairwise unions are $V(D)\times V(D)$.
Having imposed an extra condition on the minimization, we have $f^*(D)\ge f(D)$;
we prove that $l(D)\ge f^*(D)$.

On the upper side, we study the related {\it catch leafage} $l^*(D)$ of a \dig\
$D$.  This is the minimum number of leaves in a host tree in which $D$ has a
subtree \rep\ such that each sink subtree is a single vertex.  (Such \rep s,
particularly when the host tree is a path, are studied in [12,13,15].)
This condition restricts the allowable \rep s, so $l^*(D) \ge l(D)$.  We prove
that $l^*(D)\le w(P(D))$, where $w(P(D))$ is the width of the inclusion poset
$P(D)$ on the sets whose incidence vectors are the columns of the adjacency
matrix $A(D)$.  We also give a sufficient condition for equality in this bound.

We thus obtain the chain of inequalities
$$f(D)\le f^*(D)\le l(D)\le l^*(D)\le w(P(D))\le n(D) .$$
We present examples to show that each inequality is best possible.
We also present examples to show that each bound is arbitrarily weak,
as any one of these parameters can be at most 3 when the next parameter is
arbitrarily large.
%(at least $n(D)-1$ except for the first inequality).

The upper bound $w(P(D))$ is easily computable, but the lower bounds are not.
Cogis [2] and Doignon, Ducamp, and Falmagne [4] proved an easily testable
\charzn\ of the \dig s with \fd\ at most 2, but Yannakakis [18] proved that
recognition of \fd\ 3 is NP-complete.  M\"uller [11] found a polynomial-time
recognition algorithm for interval digraphs (leafage 2).  Other than this, we do
not know the complexity of recognizing digraphs with bounded values for any of
$\{f^*(D),l(D),l^*(D)\}$.

\SH
{2. SUBTREE REPRESENTATIONS AND LEAFAGE}
We use $u \to v$ to denote the successor relation; $u\to v$ means ``$uv$
is an edge''.  A {\it branch point} of a tree is a vertex of degree at least 3.
We show first that leafage is well-defined.

\TH 1.
If $D$ is a digraph with $n$ vertices, then $D$ has a subtree \rep\ in a star
with at most $n$ leaves.
\PF
In a star $H$ with $n$ leaves, assign distinct leaves as sink sets for the $n$
vertices.  For each $v\in V(D)$, let $S_v$ be the star induced by the
center of $H$ and the leaves corresponding to the successors of $v$.
Then $u\to v$ if and only if $T_v \esub S_u$, and hence this is a \rep.  \qed

The bound $l(D)\le n(D)$ is sharp, as it holds with equality for the digraph
$D_n$ in Theorem 2.  Our tool for proving lower bounds on $l(D)$ is a property
of subtrees of a tree.  If $T_i$, $T_j$, $T_k$ are subtrees of a tree, then we
say that $T_k$ is {\it between} $T_i$ and $T_j$ if $T_i\cap T_j = \nul$ and the
unique path from $T_i$ to $T_j$ contains a vertex of $T_k$ (possibly at the
start or end).  A collection of pairwise disjoint
subtrees having the property that none is between two others is an {\it
asteroidal collection} of subtrees.

\LM 1.
If $\VEC T1n$ is an asteroidal collection of subtrees of a tree $T$,
then $T$ has at least $n$ leaves.
\PF
We may assume that the path from any leaf of $T$ to the nearest branch point
contains a vertex of some $T_i$; otherwise, we could delete the vertices before
the branch point to reduce the number of leaves without changing the hypotheses.
For each leaf $v$ of $T$, we assign to $v$ the first subtree encountered on the
path from $v$ to its nearest branch point.  If $T$ has fewer than $n$ leaves,
then some subtree $T_k$ in our list is not assigned to any leaf.  Let $x$ be a
vertex of $T_k$, and let $P$ be a maximal path containing $x$.  The endpoints of
$P$ are leaves of the tree, and $T_k$ is between the subtrees assigned to those
leaves.  Hence $T$ must have at least $n$ leaves.  \qed\upsp

\LM 2.
If $v,w$ have a common successor $u$ that is not a successor of $z$ in a 
digraph $D$, then $S_z$ is not between $S_v$ and $S_w$ in any subtree \rep\ of
$D$.  Similarly, if $v,w$ have a common predecessor $u$ that is not a
predecessor of $z$ in $D$, then $T_z$ is not between $T_v$ and $T_w$ in any
subtree \rep\ of $D$.  
\PF
If $S_z$ is between $S_v$ and $S_w$, then $S_v\cap S_w = \nul$, and $T_u$
must contain the unique path from $S_v$ to $S_w$ in the host.  This
contradicts $S_z\cap T_u = \nul$, since $S_z$ has a vertex on this path.
The proof of the other statement is similar.  \qed

Subtrees of a tree satisfy the {\it Helly property}; the members of a
pairwise \intg\ family of (sub)trees have a common vertex (see, for example,
[6, p. 92]).

\LM 3.
If in a subtree \rep\ of $D$ the source subtrees are pairwise \intg\ and
the sink subtrees are pairwise \intg, then $A(D)$ has a row of 1's or a
column of 0's, and similarly $A(D)$ has a column of 1's or a row of 0's.
\PF
In such a \rep, the source subtrees have a common vertex, and the sink subtrees
have a common vertex.  Let $s,t$ denote these vertices, respectively.  If $s=t$,
then $A(D)$ is all 1's and the claim holds.  If $s \ne t$, let $x$ be the vertex
of $\bigcup S_i$ that is closest to $t$ on the unique $s,t$-path in $T$.
Suppose $x \in S_k$.  If $A(D)$ has no row of 1's, then $x\ne t$ and some sink
subtree $T_j$ fails to contain $x$.  However, $t\in T_j$, and hence $T_j$
intersects no source subtree, forcing a column of 0's in $A(D)$.  The other
claim follows by considering the vertex of $\bigcup T_j$ that is closest to $s$
on the $s,t$-path in $T$.  \qed

\medskip
Because permuting rows or columns is simply a relabeling of source or sink
sets, leafage can be viewed as a property of a 0,1-matrix (the adjacency matrix
$A(D)$) rather than a property of a \dig.  We next show that asteroidal
collections are forced by complements of permutation matrices.  

\medskip
\TH 2.
For $n\ge 3$,
let $D_n$ be the digraph such that $A(D_n)=J-I$, where $J$ is the matrix of all
1's and $I$ is the identity matrix.  In every subtree \rep\ of $D_n$, either the
source subtrees have a common vertex and the sink subtrees form an asteroidal
collection, or the sink subtrees have a common vertex and the source subtrees
form an asteroidal collection. 
\PF
Let the vertices of $D_n$ be $\{\LIST 1n\}$; we have $i\to j$ if
and only if $i\ne j$.  Consider a subtree \rep\ of $D_n$.  By Lemma 3, the
source subtrees and sink subtrees cannot both be pairwise \intg; we may assume
by symmetry that there is a disjoint pair of source subtrees.

When $i,j,k$ are distinct vertices, we have $i\to k$, $j\to k$, and
$k\not\to k$.  Thus Lemma 2 implies that no source subtree can be between two
other source subtrees.  With betweenness forbidden, $S_i$ and $S_k$ cannot 
intersect if $S_i\cap S_j=\nul$.  We conclude that if some pair of source
subtrees is disjoint, then the source subtrees are pairwise disjoint, and
none is between two others.  Hence they form an asteroidal collection.

With the source subtrees pairwise disjoint, consider the sink subtrees.  For any
distinct vertices $i,j,k$, we must have $T_j$ containing the path from $S_i$ to
$S_k$ and $T_i$ containing the path from $S_j$ to $S_k$.  Hence
$T_i\cap T_j \ne\nul$, and the sink subtrees are pairwise \intg.  The Helly
property then implies that the sink subtrees have a common vertex.  \qed

\medskip
Together, Lemma 1 and Theorem 2 imply that $l(D_n) = n$.

\SH
{3. LEAFAGE AND DISJOINT FERRERS DIMENSION}
We next prove our main lower bound on leafage.  We use $N_D^+(u)$ to denote
the successor set and $N_D^-(u)$ to denote the predecessor set of a vertex $u$
in a digraph $D$.

\TH 3.
If $D$ is a \dig, then $l(D)\ge f^*(D)$.
\PF
Suppose that $l(D)=k$, and let $\{(S_v,T_v)\st v\in V(D)\}$ be a \rep\ of $D$ in
a host tree with $k$ leaves.  When $k=2$, the result follows from the \charzn\
of \idig s in [15].  For $k\ge 3$, we construct $k$ pairwise disjoint \fdig s
whose union is $\Db$.  With the host tree $T$ embedded in the plane, let the
leaves be $\VEC x1n$ in clockwise order around the tree.  Let $P_i$ denote the
path in $T$ from $x_i$ to $x_{i+1}$, indexed cyclically.

For each leaf $x_i$ of the host tree $T$, we construct an associated \fdig\
$D(i)$.  The edges of $\Db$ consist of those pairs $uv$ such that
$S_u\cap T_v=\nul$, meaning that the unique shortest path from $S_u$ to $T_v$
has length at least 1.  Let $D(i)$ consist of those edges $uv$ in $\Db$ such
that the first edge on the path from $S_u$ to $T_v$ lies on $P_i$, with $S_u$
between $x_i$ and $T_v$ (see Fig.~1).

If $S_u$ has no vertex on $P_i$, then $u$ has no successors in $D(i)$.
If the last vertex of $S_{u'}$ on $P_i$ is closer to $x_{i+1}$ than the
last vertex of $S_u$ on $P_i$, then $N_{D(i)}^+(u')\esub N_{D(i)}^+(u)$, by
construction.  Hence the $D(i)$'s are \fdig s.

The paths $P_i$ together cover each edge of the host tree exactly once in each
direction.  Since each edge is covered in each direction, $\bigcup_i D(i)=\Db$.
Since each edge is covered only once, and when $S_u\cap T_v = \nul$ there is a
unique first edge on the path from $S_u$ to $T_v$, the subgraphs $\{D(i)\}$ 
are pairwise disjoint.  \qed\nosp
\nosp
\gpic{
\expandafter\ifx\csname graph\endcsname\relax \csname newbox\endcsname\graph\fi
\expandafter\ifx\csname graphtemp\endcsname\relax \csname newdimen\endcsname\graphtemp\fi
\setbox\graph=\vtop{\vskip 0pt\hbox{%
    \graphtemp=.5ex\advance\graphtemp by 0.239in
    \rlap{\kern 2.894in\lower\graphtemp\hbox to 0pt{\hss $\bullet$\hss}}%
    \graphtemp=.5ex\advance\graphtemp by 1.832in
    \rlap{\kern 2.894in\lower\graphtemp\hbox to 0pt{\hss $\bullet$\hss}}%
    \graphtemp=.5ex\advance\graphtemp by 1.832in
    \rlap{\kern 0.239in\lower\graphtemp\hbox to 0pt{\hss $\bullet$\hss}}%
    \graphtemp=.5ex\advance\graphtemp by 0.106in
    \rlap{\kern 0.106in\lower\graphtemp\hbox to 0pt{\hss $\bullet$\hss}}%
    \graphtemp=.5ex\advance\graphtemp by 0.239in
    \rlap{\kern 3.000in\lower\graphtemp\hbox to 0pt{\hss $x_1$\hss}}%
    \graphtemp=.5ex\advance\graphtemp by 1.832in
    \rlap{\kern 3.000in\lower\graphtemp\hbox to 0pt{\hss $x_2$\hss}}%
    \graphtemp=.5ex\advance\graphtemp by 1.832in
    \rlap{\kern 0.133in\lower\graphtemp\hbox to 0pt{\hss $x_3$\hss}}%
    \graphtemp=.5ex\advance\graphtemp by 0.106in
    \rlap{\kern 0.000in\lower\graphtemp\hbox to 0pt{\hss $x_4$\hss}}%
    \special{pn 8}%
    \special{pa 2894 239}%
    \special{pa 2097 1035}%
    \special{pa 2894 1832}%
    \special{fp}%
    \special{pa 2097 1035}%
    \special{pa 1035 1035}%
    \special{pa 239 1832}%
    \special{pa 1035 1035}%
    \special{pa 106 106}%
    \special{fp}%
    \graphtemp=.5ex\advance\graphtemp by 1.699in
    \rlap{\kern 1.566in\lower\graphtemp\hbox to 0pt{\hss $uv\in D(2),~u'v\in D(3)$\hss}}%
    \graphtemp=.5ex\advance\graphtemp by 1.938in
    \rlap{\kern 1.566in\lower\graphtemp\hbox to 0pt{\hss $N_{D(2)}^+(u')\esub N_{D(2)}^+(u)$\hss}}%
    \graphtemp=.5ex\advance\graphtemp by 0.456in
    \rlap{\kern 2.500in\lower\graphtemp\hbox to 0pt{\hss $S_u$\hss}}%
    \graphtemp=.5ex\advance\graphtemp by 1.615in
    \rlap{\kern 0.633in\lower\graphtemp\hbox to 0pt{\hss $S_{u'}$\hss}}%
    \graphtemp=.5ex\advance\graphtemp by 0.429in
    \rlap{\kern 0.580in\lower\graphtemp\hbox to 0pt{\hss $T_v$\hss}}%
    \special{pa 2257 1301}%
    \special{pa 2044 1088}%
    \special{fp}%
    \special{pa 1088 982}%
    \special{pa 876 770}%
    \special{fp}%
    \special{pa 2257 770}%
    \special{pa 2044 982}%
    \special{fp}%
    \special{pa 1088 1088}%
    \special{pa 876 1301}%
    \special{fp}%
    \special{pa 2044 1088}%
    \special{pa 1471 1088}%
    \special{fp}%
    \special{sh 1.000}%
    \special{pa 1577 1115}%
    \special{pa 1471 1088}%
    \special{pa 1577 1062}%
    \special{pa 1577 1115}%
    \special{fp}%
    \special{pa 1471 1088}%
    \special{pa 1088 1088}%
    \special{fp}%
    \special{pa 1088 982}%
    \special{pa 1662 982}%
    \special{fp}%
    \special{sh 1.000}%
    \special{pa 1556 956}%
    \special{pa 1662 982}%
    \special{pa 1556 1009}%
    \special{pa 1556 956}%
    \special{fp}%
    \special{pa 1662 982}%
    \special{pa 2044 982}%
    \special{fp}%
    \graphtemp=.5ex\advance\graphtemp by 1.221in
    \rlap{\kern 1.566in\lower\graphtemp\hbox to 0pt{\hss $P_2$\hss}}%
    \graphtemp=.5ex\advance\graphtemp by 0.850in
    \rlap{\kern 1.566in\lower\graphtemp\hbox to 0pt{\hss $P_4$\hss}}%
    \special{pa 2363 850}%
    \special{pa 2177 1035}%
    \special{fp}%
    \special{pa 956 1035}%
    \special{pa 770 1221}%
    \special{fp}%
    \special{pa 2177 1035}%
    \special{pa 2288 1147}%
    \special{fp}%
    \special{sh 1.000}%
    \special{pa 2232 1053}%
    \special{pa 2288 1147}%
    \special{pa 2195 1091}%
    \special{pa 2232 1053}%
    \special{fp}%
    \special{pa 2288 1147}%
    \special{pa 2363 1221}%
    \special{fp}%
    \special{pa 956 1035}%
    \special{pa 844 924}%
    \special{fp}%
    \special{sh 1.000}%
    \special{pa 901 1018}%
    \special{pa 844 924}%
    \special{pa 938 980}%
    \special{pa 901 1018}%
    \special{fp}%
    \special{pa 844 924}%
    \special{pa 770 850}%
    \special{fp}%
    \graphtemp=.5ex\advance\graphtemp by 1.035in
    \rlap{\kern 2.363in\lower\graphtemp\hbox to 0pt{\hss $P_1$\hss}}%
    \graphtemp=.5ex\advance\graphtemp by 1.035in
    \rlap{\kern 0.770in\lower\graphtemp\hbox to 0pt{\hss $P_3$\hss}}%
    \special{pn 28}%
    \special{pa 770 770}%
    \special{pa 239 239}%
    \special{fp}%
    \special{pa 2097 1035}%
    \special{pa 2628 504}%
    \special{fp}%
    \special{pa 1035 1035}%
    \special{pa 504 1566}%
    \special{fp}%
    \hbox{\vrule depth1.938in width0pt height 0pt}%
    \kern 3.000in
  }%
}%
}
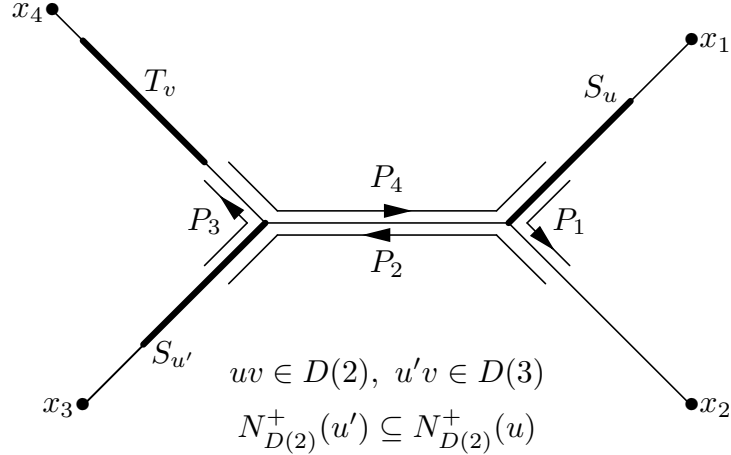
\ce{Fig.~1.  Ferrers digraphs from subtree representation}
\sp

This provides another proof that the leafage of the \dig\ $D_n$ is $n$.  Since
each pair of ones on the diagonal of $A(\Db_n)$ induce a 2 by 2 permutation
submatrix, no pair of them can be covered by a single \fdig\ contained in
$\Db_n$.

Although the inequalities $f(D)\le f^*(D)\le l(D)\le n(D)$ are best possible,
with equality throughout when $D = D_n$, the gaps can be arbitrarily large.
For an \idig, $f(D) = f^*(D) = l(D) = 2$.  By the \charzn\ of interval digraphs
in [15], $f^*(D)=2$ implies $l(D)=2$.  Nevertheless, there exist digraphs
$D$ with $f^*(D)=3$ and $l(D)=n(D)$.

\sp
\TH 4.
Leafage is not bounded by any function of $f^*$ when $f^*\ge3$.  In particular,
let $E_n$ be the $n$-vertex digraph with $A(E_n)=\CH{I\ \ Y}{Y^T\ 0}$,
%\pmatrix{I&Y\cr Y^T&0\cr},
where $I$ denotes the $n-1$ by $n-1$ identity matrix and Y denotes a column
vector of $n-1$ ones.  If $n\ge 3$, then $l(E_n) = n$, but $f^*(E_n)=f(E_n)=3$.
\PF
Because the last three rows and columns of $A(E_n)$ form a row permutation of
$A(D_3)$, we have $f^*(E_n)\ge f(E_n)\ge3$.  For equality, partition the zeros
of $A(\Eb_n)$ into three sets; those in the upper right of the submatrix $I$,
those in the lower left of the submatrix $I$, and the 0 in the lower right
corner.  These sets yield Ferrers digraphs, so $f^*(E_n)\le3$.

To show that $l(E_n) = n$, we name the vertices by the row and column indices of
the matrix and let $\{(S_i,T_i)\st 1\le i\le n\}$ be a subtree \rep\ of $E_n$
in the host tree $T$.  By Lemma 1, it suffices to show that the source
subtrees or the sink subtrees form an asteroidal collection in $T$.

We have $S_n\cap T_n = \nul$; let $P$ be the unique path from $S_n$ to $T_n$ in
$T$.  For each $k<n$, we have $S_k\cap T_n \ne \nul$, $S_n\cap T_k\ne\nul$,
and $S_k\cap T_k\ne\nul$.  Consider also $i<n$.  If $P$ contains a vertex $x$ of
$S_i\cap T_i$, then the nonadjacency of $i$ and $k$ implies that $x$ separates
$S_k$ and $T_k$.  This contradicts $S_k\cap T_k\ne\nul$, so $S_i$ cannot
intersect $T_i$ in $P$.  We conclude that $S_i\cap T_i$ is contained in the
component of $T-E(P)$ containing $T_n$ or in the component of $T-E(P)$
containing $S_n$.  By symmetry, we may assume the latter.
Since $n\to i$, we now have $P\subset T_i$.  Applying this argument for all
vertices other than $n$ yields that all $S_i\cap T_i$ lie in the same component
of $T-E(P)$, since there are no edges except loops among these vertices.
Thus $P\subset T_i$ and $P\cap S_i=\nul$ for all $i<n$. 
%Suppose some vertex $x$ of $P$ belongs to both $S_i$ and $T_j$, with $i,j<n$.
%This requires $i=j$.

Now consider disjointness and betweenness of the source subtrees.  Since
$i\to i$, $n\to i$, and $j\not\to i$, Lemma 2 forbids $S_j$ between $S_i$ and
$S_n$ for $i,j<n$.  Since $P$ separates $S_n$ from the others, this implies that
the source subtrees are pairwise disjoint.
Furthermore, if $S_j$ is between $S_i$ and $S_k$ for $i,j,k<n$, then the union
of the paths from $S_n$ to the trees $S_i$ and $S_k$ must intersect $S_j$, which
puts $S_j$ between $S_n$ and one of $\{S_i,S_k\}$.  Hence the source subtrees
are pairwise disjoint, and none is between two others.  They form an asteroidal
collection, and Lemma 1 applies.
\qed

Every $n$ by $n$ (adjacency) matrix with leafage $n$ is a minimal forbidden
submatrix for leafage less than $n$.  We next present another such family.
Given the adjacency matrix $A(D)$ of a digraph $D$, let $H(D)$ be the graph
with vertices corresponding to the zeros of $A(D)$ and edges corresponding
to the pairs of zeros contains in a 2 by 2 permutation submatrix.
Cogis [2] and Doignon-Ducamp-Falmagne [4] proved that $D$ has \fd\ 2 if and
only if $H(D)$ is bipartite; here we need only the obvious necessity of the
condition.

\TH 5.
Let $C_n$ be the digraph consisting of a directed cycle of length $n$ plus
a loop at each vertex.  Then $l(C_n) = n$, but $f(C_n)=f^*(C_n)=3$.
\PF
Assume that the cycle is $1 \to 2 \to\cdots\to n\to 1$.  Partition the zeros
of $A(C_n)$ into three sets: those in the last row, those in the first $n-1$
rows below the diagonal, and the remainder.  These sets form \fdig s,
so $f^*(C_n)\le 3$.  To prove that $f(C_n)>2$, we observe that the positions
$$\{(i,i+\CL{n/2})\st 1\le i\le\FL{n/2}\} \cup
\{(i,i+1-\CL{n/2})\st \CL{n/2}\le i\le n\}$$
form an odd cycle in $H(C_n)$.

We use induction on $n$ to prove that $l(C_n) = n$.  The claim holds for $n=3$
because $A(C_3)$ is a permutation of $A(D_3)$.  For $n>3$, let $\B T$ be the
host tree for an optimal \rep\ of $C_n$.  Suppose first that
$S_{i-1}\cap S_i \ne\nul$ for some $i$ (all indexing is circular).  The subtree
$T_i$ must intersect both of these, so by the Helly property $S_{i-1},T_i,S_i$
have a common vertex $x$ in $\B T$.  No other source subtree intersects $T_i$,
and no other sink subtree intersects $S_{i-1}$ and $S_i$; hence no other
assigned subtree contains $x$.  Every two consecutive subtrees in the list
$T_{i+1},S_{i+1},T_{i+2},\dots,S_{i-2},T_{i-1}$ intersect; hence their union is
connected and contained in one component of $\B T-x$.  The remaining components
of $\B T-x$ can be deleted without changing the intersection digraph, so we may
assume that $x$ is a leaf.

Let $P$ be the path in $\B T$ from $x$ to the nearest branch point.  By
symmetry, we may assume that $S_{i-1}$ contains as much of $P$ as $S_i$.  If
$S_i$ does not contain all of $P$, then $T_{i+1}$ intersects $S_{i-1}$, which is
forbidden.  Hence $P \esub S_i\cap S_{i-1}$, and no sink subtree other than
$T_i$ intersects $P$.  If another source subtree extends onto $P$, then deleting
its edges on $P$ does not change the \int\ digraph.  We can now delete $T_i$
and replace $S_{i-1},S_i$ by a single source subtree with edge set
$(E(S_{i-1})\cup E(S_i))-E(P)$ to obtain a \rep\ of $C_{n-1}$ with $l(C_n)-1$
leaves.  By the induction hypothesis, this yields $l(C_n)\ge n$.

Hence we may assume that $S_{i-1}\cap S_i=\nul$ for all $i$, and by symmetry
also $T_{i-1}\cap T_i=\nul$ for all $i$.  In this case, let $P_i$ be the portion
of $S_i$ that is the unique $T_i,T_{i+1}$-path, and let $Q_i$ be the portion of
$T_i$ that is the unique $S_{i-1},S_i$-path.  Note that $Q_i \cap P_i$ and
$P_i\cap Q_{i+1}$ are single vertices.  The union of all these paths is thus a
closed walk in which no consecutive edges are the same.  Such a walk contains a
cycle, which is impossible in a host tree.  Hence this case does not arise.
\qed

We have presented examples with fixed $f^*(D)$ and large $l(D)$.  Also $f^*(D)$
may be arbitrarily large when $f(D)=2$.  We construct a two-parameter family of
adjacency matrices.  The matrix $M_{k,m}$ is a $km$ by $km$ matrix consisting of
$k$ rows and $k$ columns of $m$ by $m$ blocks.  The diagonal blocks are the
identity matrix, the blocks below the diagonal consist entirely of 0's, and the
blocks above the diagonal consist entirely of 1's.  The zeros can be covered by
two Ferrers digraphs, each consisting of all the subdiagonal blocks and half of
each diagonal block; hence $f(M_{k,m})=2$.  We will prove that
$f^*(M_{k,m})=c+1$ when $k=1+\CH c2$ and $m$ is sufficiently large.
(In this discussion we use the notation $M_{k,m}$ for both the digraph and
its adjacency matrix.)

Let $I_n$ denote the $n$-vertex digraph whose adjacency matrix is the identity.
A partition of $\Ib_n$ into $c$ Ferrers digraphs can be viewed as a special
$c$-coloring of the 0's in the $n$ by $n$ identity matrix $I_n$.  We say that
colors $A,B$ are a {\it crossed pair} if $A,B$ appear together in some row and
appear together in some column.

\LM 4.
If $n\ge 3c!/2$, then every partition of $\Ib_n$ into $c$ Ferrers digraphs has a
crossed pair of colors.
\PF
The proof is by induction on $c$.  For $c=2$, a 2-coloring of the 0's in
the 3 by 3 identity matrix cannot have all rows or all columns monochromatic
without having a 2 by 2 permutation matrix with 0's in one color.  For $c>2$,
let $n=3c!/2$ and $r=3(c-1)!/2$.  Consider a partition of $\Ib_n$ into $c$
Ferrers digraphs, and suppose that the corresponding coloring has no crossed
pair.

Since each row of the identity matrix has $n-1$ 0's, the pigeonhole principle
implies that each row has at least $\CL{(3(c-1)!/2)(c/c)-1/c}=r$ 0's in some
color.  By symmetry, we may assume there are 0's of color $A$ in the first $r$
columns of row $r+1$ (see Fig.~2).  Let $D$ be the subdigraph induced by the
first $r$ vertices, with $K$ the corresponding submatrix.  By the induction
hypothesis, every partition of $\Db$ into $c-1$ Ferrers digraphs yields a
coloring of the 0's in $K$ with a crossed pair of colors.  Hence we may assume
that all $c$ colors (including $A$) appear in $K$.

Let $i$ be the index of a row in $K$ in which color $A$ appears.  If another
color appears in row $i$ of $K$, then it crosses $A$ in the full matrix.  Thus
we may assume that row $i$ of $K$ has only color $A$.  Now, to avoid the
forbidden submatrix in color $A$, position $i,r+1$ must have some other color
$B$.  Now colors $A$ and $B$ appear in a row together, so they cannot appear in
a column together.  This contradicts the observation that every color, including
$B$, appears in $K$.  \qed\nosp
\gpic{
\expandafter\ifx\csname graph\endcsname\relax \csname newbox\endcsname\graph\fi
\expandafter\ifx\csname graphtemp\endcsname\relax \csname newdimen\endcsname\graphtemp\fi
\setbox\graph=\vtop{\vskip 0pt\hbox{%
    \graphtemp=.5ex\advance\graphtemp by 1.175in
    \rlap{\kern 0.343in\lower\graphtemp\hbox to 0pt{\hss $A$\hss}}%
    \graphtemp=.5ex\advance\graphtemp by 1.175in
    \rlap{\kern 0.524in\lower\graphtemp\hbox to 0pt{\hss $A$\hss}}%
    \graphtemp=.5ex\advance\graphtemp by 1.175in
    \rlap{\kern 0.705in\lower\graphtemp\hbox to 0pt{\hss $A$\hss}}%
    \graphtemp=.5ex\advance\graphtemp by 1.175in
    \rlap{\kern 0.886in\lower\graphtemp\hbox to 0pt{\hss $A$\hss}}%
    \graphtemp=.5ex\advance\graphtemp by 1.175in
    \rlap{\kern 1.066in\lower\graphtemp\hbox to 0pt{\hss $A$\hss}}%
    \graphtemp=.5ex\advance\graphtemp by 1.175in
    \rlap{\kern 1.247in\lower\graphtemp\hbox to 0pt{\hss $A$\hss}}%
    \graphtemp=.5ex\advance\graphtemp by 0.090in
    \rlap{\kern 0.343in\lower\graphtemp\hbox to 0pt{\hss $1$\hss}}%
    \graphtemp=.5ex\advance\graphtemp by 0.271in
    \rlap{\kern 0.524in\lower\graphtemp\hbox to 0pt{\hss $1$\hss}}%
    \graphtemp=.5ex\advance\graphtemp by 0.452in
    \rlap{\kern 0.705in\lower\graphtemp\hbox to 0pt{\hss $1$\hss}}%
    \graphtemp=.5ex\advance\graphtemp by 0.633in
    \rlap{\kern 0.886in\lower\graphtemp\hbox to 0pt{\hss $1$\hss}}%
    \graphtemp=.5ex\advance\graphtemp by 0.813in
    \rlap{\kern 1.066in\lower\graphtemp\hbox to 0pt{\hss $1$\hss}}%
    \graphtemp=.5ex\advance\graphtemp by 0.994in
    \rlap{\kern 1.247in\lower\graphtemp\hbox to 0pt{\hss $1$\hss}}%
    \graphtemp=.5ex\advance\graphtemp by 0.452in
    \rlap{\kern 0.343in\lower\graphtemp\hbox to 0pt{\hss $A$\hss}}%
    \graphtemp=.5ex\advance\graphtemp by 0.452in
    \rlap{\kern 0.524in\lower\graphtemp\hbox to 0pt{\hss $A$\hss}}%
    \graphtemp=.5ex\advance\graphtemp by 1.175in
    \rlap{\kern 1.428in\lower\graphtemp\hbox to 0pt{\hss $1$\hss}}%
    \graphtemp=.5ex\advance\graphtemp by 0.452in
    \rlap{\kern 0.886in\lower\graphtemp\hbox to 0pt{\hss $A$\hss}}%
    \graphtemp=.5ex\advance\graphtemp by 0.452in
    \rlap{\kern 1.066in\lower\graphtemp\hbox to 0pt{\hss $A$\hss}}%
    \graphtemp=.5ex\advance\graphtemp by 0.452in
    \rlap{\kern 1.247in\lower\graphtemp\hbox to 0pt{\hss $A$\hss}}%
    \graphtemp=.5ex\advance\graphtemp by 0.452in
    \rlap{\kern 1.428in\lower\graphtemp\hbox to 0pt{\hss $B$\hss}}%
    \graphtemp=.5ex\advance\graphtemp by 0.633in
    \rlap{\kern 0.072in\lower\graphtemp\hbox to 0pt{\hss $K$\hss}}%
    \special{pn 8}%
    \special{pa 253 1084}%
    \special{pa 1337 1084}%
    \special{pa 1337 0}%
    \special{pa 253 0}%
    \special{pa 253 1084}%
    \special{fp}%
    \hbox{\vrule depth1.247in width0pt height 0pt}%
    \kern 1.500in
  }%
}%
}
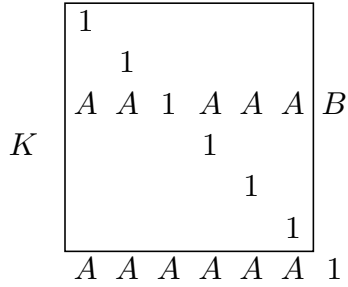
\ce{Fig.~2. Coloring 0's in an identity matrix.}
\sp

The bound $3c!/2$ in Lemma 4 is not best possible.  For $c=2,3,4$, the bound is
$3,9,36$, but the actual minimum values forcing the desired behavior are
$3,4,6$.  We are content with the bound arising from the short argument in Lemma
4 because our aim is to show that $f^*(M_{k,m})$ grows arbitrarily large.

\TH 6.
If $k\ge1+\CH c2$ and $m\ge3c!/2$, then $f^*(M_{k,m})>c$.
\PF
Suppose $\Mb_{k,m}$ has a partition into $c$ pairwise-disjoint Ferrers
digraphs.  By Lemma 4, in each copy of $I_m$ in the block structure of
$M_{k,m}$, the corresponding coloring has a crossed pair of colors.
Since there are more than $\CH c2$ diagonal blocks, by the pigeonhole
principle some pair of colors $A,B$ is crossed twice.

Let $r,s$ be the indices of the diagonal blocks where $A,B$ are crossed,
with $r<s$.  Let $j$ be the column within diagonal block $r$ where $A,B$
both appear, and let $i$ be the row within diagonal block $s$ where $A,B$
both appear.  Position $i,j$ of block $s,r$ is now forced to have both
color $A$ and color $B$ to avoid the forbidden substructure for the Ferrers
digraphs given by colors $A$ and $B$.  This is impossible.  \qed

It is worth noting that $f^*(M_{k,m})\le c$ for all $m$ when $k\le \CH c2$.
This is illustrated by the block coloring in Fig.~3.
\nosp
$$\pmatrix{A\bs B &  1  &  1  &  1  &  1  &  1  &  1  &  1  &  1  &  1\cr
            A   &A\bs C &  1  &  1  &  1  &  1  &  1  &  1  &  1  &  1\cr
            B   &  C  &B\bs C &  1  &  1  &  1  &  1  &  1  &  1  &  1\cr
            A   &  A  &  A  &A\bs D &  1  &  1  &  1  &  1  &  1  &  1\cr
            B   &  B  &  B  &  D  &B\bs D &  1  &  1  &  1  &  1  &  1\cr
            C   &  C  &  C  &  D  &  D  &C\bs D &  1  &  1  &  1  &  1\cr
            A   &  A  &  A  &  A  &  A  &  A  &A\bs E &  1  &  1  &  1\cr
            B   &  B  &  B  &  B  &  B  &  B  &  E  &B\bs E &  1  &  1\cr
            C   &  C  &  C  &  C  &  C  &  C  &  E  &  E  &C\bs E &  1\cr
            D   &  D  &  D  &  D  &  D  &  D  &  E  &  E  &  E  &D\bs E\cr}$$
\ce{Fig.~3.  A 5-coloring of the 0's in the blocks of $M_{10,l}$.}
\medskip

We previously gave examples with $f(D)=3$, $f^*(D)=3$, and $l(D)$ large.
We next prove that the family $M_{k,m}$ includes examples with $f(D)=2$,
$f^*(D)=3$, and $l(D)$ large.

\TH 7.
For $m\ge 3$, $M_{2,m}$ is a $2m$-vertex digraph with $f(M_{2,m})=2$,
$f^*(M_{2,m})=3$, and $l(M_{2,m})=m$.
\PF
With four blocks of order $m$, $M_{2,m}=\CH{I\ 1}{0\ I}$.  The value
$f(M_{2,m})=2$ was obtained before Lemma 4.  The lower bound on $f^*$ comes from
Theorem 6 (with $c=2$), and the upper bound comes from the coloring illustrated
in Fig.~3 (with $c=3$).

To prove that $l(M_{2,m})\le m$, we construct a representation with $m$ leaves.
Let the host tree be the union of $m$ paths $q,s_i,t_i$ with $q$ as a common
endpoint.  Let $S_i$ and $T_{i+m}$ be the entire $i$th path, for $1\le i\le m$.
Let $S_{i+m}$ be the single vertex $s_i$, and let $T_i$ be the single vertex
$t_i$.

It remains to prove that $l(M_{2,m})\ge m$.  Consider a subtree representation
with source subtrees $\VEC S1{2m}$ and sink subtrees $\VEC T1{2m}$ for the
vertices indexed by the rows and columns of $M_{2,m}$ in order.  If
$\VEC T{m+1}{2m}$ have no common point, then some $T_i,T_j$ among these
are disjoint.  Since $T_i$ and $T_j$ must intersect each of $\VEC S1m$, those
subtrees contain the $T_i,T_j$-path and hence have a common point.  Similarly,
if $\VEC S1m$ have no common point, then $\VEC T{m+1}{2m}$ must.  By symmetry,
we may assume that $\VEC S1m$ have a common point $q$.

We now show that $\VEC T1m$ is an asteroidal collection of subtrees.  If
$T_i\cap T_j\ne\nul$ with $i,j\le m$, then the entire path from $q$ to
the closest vertex of $T_i\cap T_j$ belongs to at least one of $\{S_i,S_j\}$,
which contradicts the requirement that each of $\{S_i,S_j\}$ intersects exactly
one of $\{T_i,T_j\}$.  If $T_j$ is between $T_i$ and $T_k$, then let $P$ be the
path between $T_i$ and $T_k$, and let $r$ be the vertex of $P$ closest to $q$.
Depending on the location of $r$ relative to $T_j$ on $P$, the $q,T_k$-path in
$S_k$ or the $q,T_i$-path in $S_i$ intersects $T_j$, contradicting their
disjointness from $T_j$.  Thus $\VEC T1m$ is an asteroidal collection,
and Lemma 1 implies that the host has at least $m$ leaves.  \qed

\SH
{4. CATCH LEAFAGE}
If $D$ has a subtree \rep\ in which every sink subtree is a single vertex,
then we say this is a {\it catch \rep}, and $D$ is a {\it catch-tree digraph}.
In discussing catch representations, we say ``sink point'' instead of ``sink
subtree'' to make the usage clear.  If $D$ has a catch-tree \rep\ in which the
host is a path, then $D$ is a {\it catch-interval \dig}.  The corresponding
classes in which the source sets are single vertices are merely those whose
adjacency matrices are the transposes of the digraphs in the classes defined
above.  Catch-interval \dig s are characterized in [12] under the name
``interval catch digraphs'' and in [15] under the name ``interval-point
\dig s''.

The {\it catch leafage} $l^*(D)$ is the minimum number of leaves in a host tree
in which $D$ has a catch-tree \rep; the catch-interval \dig s are 
the \dig s with catch-leafage 2.  In the proof of Theorem 1, we gave every
$n$-vertex \dig\ a catch \rep\ in a star with $n$ leaves, so catch leafage is
well-defined.  Since every catch-tree \rep\ is a subtree \rep, we have
$n\ge l^*(D)\ge l(D)$.

We may make several simplifying assumptions about the
form of optimal catch-tree \rep s.  In a catch-tree \rep, sink subtrees can
occupy the same vertex of the host if and only if the corresponding columns of
the matrix are identical.  We may split such a vertex of the host (without
increasing the number of leaves), including the source subtrees to cover both.
Thus we may assume that in catch representations each vertex is occupied by at
most one sink point.  Also, if a vertex of degree at most two in the host tree
is not assigned as a sink point, then an edge incident to it can be contracted.

Recall that the predecessor set for $v$ is $N^-(v)=\{u\st u\to v\}$; this is the
set whose incidence vector is the column of the adjacency matrix
corresponding to $v$.  Because the source sets occupy single vertices, a
catch-tree \rep\ can be described by listing, for each vertex of the host tree,
the non-empty collection of source sets containing it.  This will be a
catch-tree \rep\ if and only if 1) among these sets appear the predecessor sets,
and 2) the set of vertices assigned to each source set forms a subtree of the
host.

Therefore, our analysis of catch leafage focuses on the columns of the adjacency
matrix as incidence vectors for the predecessor sets.  We define an associated
partial order.  Let $P(D)$, the {\it incidence poset} of the digraph $D$, be the
collection of predecessor sets in $D$, ordered by inclusion.  For simplicity, we
will use the same notation $V_j$ to refer to a predecessor set or the column of
the adjacency matrix that is its incidence vector.
%and we call the elements of the predecessor sets (i.e., the rows of the
%adjacency matrix) the {\it predecessors}.

The {\it width} $w(P)$ of a poset $P$ is the maximum size of its antichains
(collections of pairwise incomparable elements).  Dilworth's Theorem [3] says
that the elements of a finite poset $P$ can be partitioned into $w(P)$
disjoint chains.

\TH 8.
The inequality $l^*(D)\le w(P(D))$ holds for every \dig\ $D$ with
$w(P(D))\ge 2$.
\PF
Let $k=w(P(D)$, and let $\VEC C1k$ be a partition of $P(D)$ into $k$ disjoint
chains.  Let the host tree $T$ be a subdivision of a star with $k$ leaves.
That is, $T$ consists of a central point of degree $k$ from which $k$ paths
emerge.  Assign the central vertex the set of all predecessors, and
assign to each emerging path the sets on one of the chains $C_i$, in
decreasing order.  The predecessor sets all appear at vertices, and the
occurrences of each predecessor form a subtree, so this is a catch-tree
\rep.  \qed

Fulkerson [5] observed that Dilworth's Theorem is equivalent to the
K\"onig-Egerv\'ary Theorem on matchings in bipartite graphs.  Thus bipartite
matching or other algorithms can be used to compute $w(P(D))$.  Nevertheless,
this is only a bound on $l^*(D)$, and this bound also can be arbitrarily bad.
The digraph $D$ consisting of a directed path plus a loop at each vertex has
catch leafage 2 but $w(P(D)) = n-1$, so $w(P(D))$ is not bounded by any
function of $l^*(D)$.

Note that $w(P(D))=1$ when $D$ is a Ferrers digraph.  Thus Theorem 8 requires
$w(P(D))\ge 2$, and we see that the break between $w(P(D))$ and $n(D)$ can
be large.

We now have the chain of inequalities 
\ce{
$f(D)\le f^*(D)\le l(D)\le l^*(D)\le w(P(D)\le n(D)$.
}
\LP
One may have equality throughout (achieved by $D_n$).  To prove that there can
be arbitrarily bad breaks between any pair, it suffices to produce examples
where $l(D)$ is bounded and $l^*(D)$ is large.  To do this, we prove a
sufficient condition for $l^*(D)=w(P(D))$.

\TH 9.
If $D$ is a \dig\ such that $P(D)$ has a unique maximal element and
$w(P(D))\ge2$, then $l^*(D) = w(P(D))$.
\PF
Let $V_0$ be the unique maximal element, and let $A= \VEC V1k$ denote a maximum
antichain in $P(D)$.  Let $q_i$ denote the vertex of the host assigned to
$V_i$ in an optimal catch representation.  Iteratively delete leaves of the host
tree that are not in $\{q_i\}$ until all remaining leaves belong to $\{q_i\}$.
If the number of leaves (other than $q_0$) is less than $k$, then some set $V_i$
in $A$ is assigned to a non-leaf $q_i$.  Every path from $q_0$ to another
remaining vertex can be extended to reach a remaining leaf.  In particular, the
path from $q_0$ to $q_i$ belongs to a path from $q_0$ to a leaf assigned $q_j$.
Since $V_j \esub V_0$ and each predecessor is assigned to the vertices of a
tree, this entire path including $q_i$ belongs to the source subtrees for $V_j$.
This yields $V_j\esub V_i$, contradicting the choice of $A$ as an antichain.
\qed\nosp

\TH{10}.
Catch leafage is not bounded by any function of leafage.  If $F_n$ denotes the
$n$-vertex digraph whose adjacency matrix is $\CH{I\ \ Y}{Y^T\ 1}$, where $I$
denotes the $n-1$ by $n-1$ identity matrix and Y denotes a column vector of
$n-1$ ones, then $f^*(F_n)=l(F_n)=2$, but $l^*(F_n)=n-1$.
\PF
The upper left and lower right zeros in the portion $I$ of the adjacency 
matrix yield two disjoint \fdig s whose union is $\Fb_n$.  As proved in [15],
this is equivalent to leafage 2.  On the other hand, the predecessor set of
the last vertex contains all the other predecessor sets, so
$l^*(F_n)=w(P(F_n))=n-1$.  \qed

The sufficient condition in Theorem 9 does not characterize equality in
$l^*(D)\le w(P(D))$.  For the digraph $C_n$ consisting of a directed cycle plus
loops, we have seen that $l(C_n) = n$.  Also the columns of $A(C_n)$ form an
antichain, so $l(C_n)= l^*(C_n) = w(P(C_n)) = n$.

This example shows also that leafage and catch leafage can drop arbitrarily much
when a single vertex is deleted.  Deleting one vertex from a cycle with loops
leaves a path with loops.  The former has leafage and catch leafage $n$; the
latter has leafage and catch leafage 2.

Our proof of $l^*(D)\le w(P(D))$ shows that every digraph has a catch
representation in a host tree having only one branch point, and if $P(D)$
has a unique maximum this can be achieved in a host tree with the minimum
number of leaves.  This is not true of all digraphs.  The digraph $D$
with adjacency matrix below contains $C_4$ and thus has catch leafage
at least 4.  However, every catch representation of $D$ in a host tree with four
leaves has two branch points.  We thus close by mentioning two further
optimization problems for digraphs with catch leafage $k$:  Among
catch representations in trees with $k$ leaves, what is the minimum number
of branch points, and what is the minimum number of vertices?
$$
\pmatrix{1&0&0&1&0\cr
         1&1&0&0&1\cr
         0&1&1&0&1\cr
         0&0&1&1&1\cr
         1&1&1&1&1\cr}
$$

\SH
{\ce{References}}
\frenchspacing
\BP[1]
L.W.~Beineke and C.M.~Zamfirescu, Connection digraphs and second order line
graphs.  \DM\ 39(1982), 237--254.
\BP[2]
O. Cogis, A characterization of digraphs with Ferrers dimension 2.
CNRS Research Report 19(1979).
\BP[3]
R.P. Dilworth, A decomposition theorem for partially ordered sets,
{\it Ann. Math.} 51(1950), 161--165.
\BP[4]
J.P. Doignon, A. Ducamp, and J.-C. Falmagne, On realizable biorders and the
biorder dimension of a relation.  {\it J.\ Math.\ Psych.\ } 28(1984), 73-109.
\BP[5]
D.R.~Fulkerson, Note on Dilworth's decomposition theorem for partially ordered
sets. {\it Proc. Amer.  Math. Soc.} 7(1956), 701--702.
\BP[6]
M.C.~Golumbic, {\it Algorithmic graph theory and perfect graphs},
(Academic Press, 1980).
\BP[7]
F.~Harary, J.A.~Kabell, and F.R.~McMorris, Bipartite intersection graphs.
{\it Comm.\ Math.\ Univ.\ Carolinae} 23(1982), 739--745.
\BP[8]
I.-J.~Lin, T.A.~McKee, and D.B.~West, Leafage of chordal graphs, submitted.
\BP[9]
I.-J.~Lin, M.K.~Sen, and D.B.~West, Classes of interval digraphs and 
0,1-matrices, \CNum\ 123-128(1997), to appear.
\BP[10]
I.-J.~Lin and D.B.~West, Interval digraphs that are indifference digraphs.
{\it Graph theory, Combinatorics, and Algorithms}
(Y. Alavi and A. Schwenk, eds.),  {\it Proc.\ 7th Intl.\ Conf.\ Graph Th.\ -
Kalamazoo 1992}, (Wiley 1995), 751--765.
\BP[11]
H.~M\"uller, Recognizing interval digraphs and bi-interval graphs in
polynomial time, to appear.
\BP[12]
E.~Prisner, A characterization of interval catch digraphs, {\it Discrete Math.}
73(1989), 285--289.
\BP[13]
E.~Prisner, Algorithms for interval catch digraphs, {\it 2nd Twente Workshop on
Graphs and Combinatorial Optimization (Enschede, 1991)}, {\it Discrete Appl.
Math.} 51(1994), 147--157.
\BP[14]
J.~Riguet, Les relations de Ferrers.  {\it Comptes Rendus des S\'eances
hebdomadaires de l'Acad\'emie des Sciences} (Paris) 232(1951), 1729--1730.
\BP[15]
M.~Sen, S.~Das, A.B.~Roy, and D.B.~West, Interval digraphs: an analogue of
interval graphs.  \JGT\ 13(1989), 189--202.
\BP[16]
M.~Sen, S.~Das, and D.B.~West, Circular-arc digraphs: a characterization.
\JGT\ 13(1989), 581--592.
\BP[17]
M.~Sen and B.K.~Sanyal, Indifference digraphs: a generalization of indifference
graphs.  \SIDM\ 7(1994), 157--165.
\BP[18]
M.~Yannakakis, The complexity of the partial order dimension problem,
\SIAD\ 3(1982), 351-328.

\end